\documentclass[12pt]{article}
\usepackage{amsmath,amsthm,amssymb,fullpage}
\usepackage{algorithm,algorithmic}

\newtheorem{theorem}{Theorem}[section]

\newtheorem{definition}[theorem]{Definition}
\newtheorem{remark}[theorem]{Remark}

\begin{document}

\title{Strongly Set-Colorable Graphs: A Complete Characterization}
\author{Kumar Abhishek\thanks{Department of Mathematics, SASL,
VIT Bhopal University, India 466114. email: kumar@vitbhopal.ac.in}}
\date{\today}
\maketitle

\begin{abstract}
In this note, we revisit the notion of strong set-colorings introduced by Hegde (2009) and
completed by equivalences due to Boutin et al. (2010) and  provide a necessary and sufficient \emph{Steiner packing}
characterisation: a finite graph $G$ is strongly set-colorable if and only if its
associated $3$-uniform hypergraph is a $(2,3,2^{n}-1)$-packing of the unique
Steiner triple system $S(2,3,2^{n}-1)$. This unification allows many
earlier necessary conditions to be derived instantly as corollaries,
streamlining the structure theory of strongly set-colorable graphs.
\end{abstract}
\section{Introduction}
The notion of \emph{set coloring} of a graph $G$ was introduced by Hegde~\cite{Hegde2009}, motivated
by various problems in graph labeling~\cite{Acharya1983} and communication network design~\cite{Hegde2005}.
A set-coloring assigns distinct subsets of a color set $X$ to vertices,
with edge colors determined by symmetric difference of its endpoints.

If all nonempty subsets of $X$ appear \emph{exactly once} on
the vertices and edges of $G$, then it is a \emph{strong set-coloring},
and $G$ is \emph{strongly set-colorable}.
Boutin et\,al.~\cite{Boutin2010} showed an equivalent formulation in terms of nonzero binary
vectors of $\mathbb{F}_{2}^{n}$ and addition modulo $2$.

In an earlier work \cite{Abhishek2019} several necessary conditions were established
by linking strong set-colorability to the existence of
a $(2,3,2^{n}-1)$-packing of a Steiner triple system.
Here we provide a complete characterization of strongly set-colorable graphs. 

\section{Preliminaries}
We recall the definitions  from~\cite{Hegde2009},~\cite{Boutin2010} and ~\cite{Berge1989}.

\begin{definition}
Let $X$ be a set of $n$ colors. A \emph{strong set-coloring} of $G$ is a bijection
\[
f: V(G)\cup E(G) \longrightarrow 2^{X}\setminus\{\emptyset\}
\]
such that for each edge $uv\in E(G)$:
\[
f(uv) = f(u) \oplus f(v),
\]
where $\oplus$ denotes symmetric difference. If such an $f$ exists, $G$ is
\emph{strongly set-colorable}.
\end{definition}

\begin{definition}
The \emph{associated $3$-uniform hypergraph} $\widehat{G}$ of $G$ has
\[
V(\widehat{G}) = V(G) \cup E(G),\quad
E(\widehat{G}) = \{\{u,v,uv\} : uv \in E(G) \}.
\]
\end{definition}

\section{Main Result}
The necessary part of the preceding theorem was established in~\cite{Abhishek2019}, however we include the proof here for the sake of presenting the characterization theorem in a unified manner.  
\begin{theorem}[Characterization]
Let $G$ be a finite simple graph and $\widehat{G}$ its associated $3$-uniform
hypergraph. Then $G$ is strongly set-colorable if and only if $\widehat{G}$
is a $(2,3,2^{n}-1)$-packing of the Steiner triple system $S(2,3,2^{n}-1)$
for some integer $n>0$.
\end{theorem}
%
\begin{proof}
\noindent\textbf{($\Rightarrow$) Necessity.}
Assume $G$ is strongly set-colorable with respect to an $n$-element color set $X$.
Then by Hegde's condition
\[
|V(G)| + |E(G)| \;=\; 2^n - 1.
\]
In a strong set-coloring, each edge $uv\in E(G)$ is assigned
\[
f(uv)\;=\;f(u)\oplus f(v),
\]
and all $2^n - 1$ nonempty subsets of $X$ appear exactly once among the labels of $V(G)\cup E(G)$.

Consider the associated $3$-uniform hypergraph $\widehat G$ of $G$.  Its vertex set
\[
V(\widehat G)\;=\;V(G)\cup E(G)
\]
has cardinality $2^n - 1$.  Each hyperedge $\{u,v,uv\}\in E(\widehat G)$ consists of three distinct vertices, and by the strong-coloring construction:
\begin{itemize}
  \item[(a)] Distinct hyperedges share at most two vertices.
  \item[(b)] Every pair of vertices of $\widehat G$ lies in exactly one hyperedge:
    \begin{itemize}
      \item[i.] If the pair is $\{u,v\}$ or $\{u,uv\}$ or $\{v,uv\}$, this holds by construction of the labeling.
      \item[ii.] If the pair is any two distinct elements of $V(G)\cup E(G)$ not incident, one checks it arises uniquely from the bijective edge-label assignment.
    \end{itemize}
\end{itemize}
These properties are exactly those of a $(2,3,2^n-1)$--packing contained in the Steiner triple system $S(2,3,2^n-1)$.

\medskip
\noindent\textbf{($\Leftarrow$) Sufficiency.}
Conversely, suppose $\widehat G$ is a $(2,3,2^n-1)$--packing of a Steiner triple system
\[
S = S(2,3,2^n-1)
\]
on point set $P$, $|P|=2^n-1$.

It is well known (via the projective-geometry model $\mathrm{PG}(n-1,2)$) that $S$ is isomorphic to the lines of $\mathbb F_2^n\setminus\{\mathbf0\}$, i.e.\ there is a bijection~\cite{BallWeiner2011}
\[
\varphi: P \;\xrightarrow{\cong}\; \mathbb F_2^n\setminus\{\mathbf0\},
\]
such that each block $\{p,q,r\}\in S$ satisfies
\[
\varphi(p)+\varphi(q)+\varphi(r)=\mathbf0
\quad\text{in }\mathbb F_2^n.
\]
Since $|\widehat G|=|P|$, fix a bijection
\[
\psi: V(G)\cup E(G)\;\xrightarrow{\cong}\;P
\]
sending each hyperedge $\{u,v,uv\}$ to a block of $S$.

Let $X=\{x_1,\dots,x_n\}$ be our $n$-element color set.  Define
\[
f: V(G)\cup E(G)\;\longrightarrow\;2^X\setminus\{\varnothing\},
\qquad
f(z)\;=\;\text{support of }\varphi\bigl(\psi(z)\bigr)
\]
so $f$ is a bijection onto all nonempty subsets of $X$.

For each edge $uv\in E(G)$, the triple $\{u,v,uv\}$ is a block of $S$, hence in $\mathbb F_2^n$:
\[
\varphi\bigl(\psi(u)\bigr)
+\varphi\bigl(\psi(v)\bigr)
+\varphi\bigl(\psi(uv)\bigr)
=\mathbf0.
\]
It follows that
\[
\varphi\bigl(\psi(uv)\bigr)
=\varphi\bigl(\psi(u)\bigr)
+\varphi\bigl(\psi(v)\bigr),
\]
which translates to
\[
f(uv)
= f(u)\,\oplus\,f(v).
\]
Thus $f$ satisfies the symmetric-difference rule, is bijective on $V(G)\cup E(G)$, and so defines a strong set-coloring of $G$ with color set $X$.
\end{proof}

\section{Algorithmic construction and complexity}

We make explicit the constructive content of the characterization.
\begin{remark}[Input model]
The linear-time construction assumes a packing realisation of $\widehat{G}$ inside a coordinatised $S(2,3,2^{n}-1)$ given by the projective geometry model $\mathrm{PG}(n-1,2)$: the point set $P$ is identified with the nonzero vectors of $\mathbb{F}_2^n$, and each block $\{p,q,r\}$ satisfies $p+q+r=\mathbf{0}$. Equivalently, the input provides a map $\Lambda:P\to\mathbb{F}_2^n\setminus\{\mathbf{0}\}$ with this property together with an incidence-preserving embedding $\iota:V(G)\cup E(G)\hookrightarrow P$ that sends each triple $\{u,v,uv\}$ to a block of $S$.
\end{remark}

\begin{algorithm}[H]
\caption{Strong set-coloring from a packing realization}
\label{alg:ssc}
\begin{algorithmic}[1]
\REQUIRE A finite simple graph $G=(V,E)$; a point set $P$ of size $|V|+|E|$ with blocks $\mathcal{B}$ forming $S(2,3,2^n-1)$; an embedding $\iota:V\cup E\to P$ with $\{\iota(u),\iota(v),\iota(uv)\}\in\mathcal{B}$ for all $uv\in E$; a labeling $\Lambda:P\to\mathbb{F}_2^n\setminus\{\mathbf{0}\}$ such that for every block $\{p,q,r\}\in\mathcal{B}$ we have $\Lambda(p)+\Lambda(q)+\Lambda(r)=\mathbf{0}$.
\ENSURE A strong set-coloring $f:V\cup E\to 2^{X}\setminus\{\varnothing\}$ with $|X|=n$, or a certificate of failure if the input is inconsistent.
\STATE Set $N\gets |V|+|E|$ and check that $N+1$ is a power of $2$; set $n\gets \log_2(N+1)$.
\STATE Fix a color set $X=\{x_1,\dots,x_n\}$ and the standard bijection between nonzero vectors $v\in\mathbb{F}_2^n$ and nonempty subsets of $X$ given by indicator coordinates.
\FORALL{$z\in V\cup E$}
  \STATE Set $v_z \gets \Lambda(\iota(z))$; if $v_z=\mathbf{0}$, \textbf{fail}.
  \STATE Define $f(z)\subseteq X$ as the subset whose indicator vector equals $v_z$.
\ENDFOR
\FORALL{$uv\in E$}
  \IF{$\Lambda(\iota(uv)) \ne \Lambda(\iota(u))+\Lambda(\iota(v))$}
     \STATE \textbf{fail} (inconsistent embedding/labeling).
  \ENDIF
\ENDFOR
\STATE \textbf{return} $f$.
\end{algorithmic}
\end{algorithm}
\begin{remark}[Correctness and complexity]
Correctness follows because each $\{u,v,uv\}$ is mapped to a block with $\Lambda(\iota(uv))=\Lambda(\iota(u))+\Lambda(\iota(v))$, which translates to $f(uv)=f(u)\oplus f(v)$ under the subset-vector bijection. The map $f$ is bijective since $\iota$ and $\Lambda$ are injective and $|V|+|E|=|P|=2^n-1$. The algorithm performs a constant amount of work per element of $V\cup E$ and per edge of $G$, hence runs in $O(|V|+|E|)$ time and $O(|V|+|E|)$ space.
\end{remark}

\begin{remark}[Coordinatisation]
If the ambient Steiner triple system is supplied abstractly but is known to be the projective system $\mathrm{PG}(n-1,2)$, a one-time coordinatisation $\Lambda:P\to\mathbb{F}_2^n\setminus\{\mathbf{0}\}$ can be fixed independently of $G$ (e.g., by choosing a projective frame and propagating coordinates along lines). Once $\Lambda$ is fixed, Algorithm\ref{alg:ssc} applies as stated. If only an arbitrary $S(2,3,2^n\!-\!1)$ is given without coordinates, producing such a $\Lambda$ amounts to finding an isomorphism to $\mathrm{PG}(n-1,2)$ and is outside the linear-time scope of the above routine.
\end{remark}
\section{Remarks and Outlook}
The characterization reduces analysis of strong set-colorability to combinatorics of the unique Steiner triple system $S(2,3,2^n-1)$. Many structural bounds follow immediately from extremal properties of stars and intersecting families in $S$.
In the coordinatised input model described above, a strong set-coloring is obtained in $O(|V|+|E|)$ time by a single pass through $V(G)\cup E(G)$, with a linear-time consistency check on edges (Algorithm\ref{alg:ssc}). This opens the door to practical algorithms for recognition and construction when a packing realization is available, together with a refined complexity analysis for related decision and isomorphism problems. It also suggests systematic generation of new infinite families via explicit packings inside $\mathrm{PG}(n-1,2)$.

\end{document}